\documentclass[a4,11pt]{article}
\usepackage{amssymb,amsfonts,amsthm,amsmath}
\usepackage[brazil]{babel}
\usepackage{epsfig}
\usepackage[utf8]{inputenc}
\usepackage{mathrsfs}
\usepackage{dsfont}
\newtheorem{definition}{Definition}

\newtheorem{application}{Application}
\newtheorem{property}{Property}

\usepackage{graphicx}
\usepackage{amssymb}
\usepackage{pgfplots}
\usepackage{float}

\begin{document}

\title{ \Large{\bf Analytical solutions for Navier-Stokes equations with $\psi$-Caputo fractional derivative}}


\author{ {\bf {\large{D. S. Oliveira }}} \\
 {\small Coordination of Civil Engineering,} \\
{\small UTFPR,} \\
{\small 85053-525, Guarapuava, PR, Brazil}\\
 {\small oliveiradaniela@utfpr.edu.br} \\ \\
 \and
  {\bf {\large{E. Capelas de Oliveira}}} \\
 {\small Department of Applied Mathematics,} \\
{\small Imecc - Unicamp,} \\
 {\small 13083-859, Campinas, SP, Brazil} \\
 {\small capelas@unicamp.br} }

 \date{}

\maketitle

\vspace{-.4cm}


\noindent{\bf Abstract:} This work aims to use the homotopy analysis method to obtain analytical solutions of linear
time-fractional Navier-Stokes equations with cylindrical coordinates and of a system of 
nonlinear time-fractional Navier-Stokes equations with Cartesian coordinates.
These equations are described in the $\psi$-Caputo time-fractional derivative.
The solutions obtained for time-fractional Navier-Stokes equations are 
graphically presented.\\


\vspace*{.2cm}
\section{Introduction} 

Many analytical methods have been developed to solve nonlinear ordinary/par\-tial differential equations
and nonlinear fractional ordinary/partial differential equations. The method proposed by Adomian the so-called 
Adomian decomposition method (ADM) \cite{Adomian} was used by Jafari and Daftardar-Gejji \cite{Jafari} to solve 
linear and nonlinear fractional diffusion and wave equations. The homotopy perturbation method (HPM) was
constructed by He \cite{He} and was applied, recently, by Kashkari et al. \cite{Kashkari} to study 
dissipative nonplanar solitons in an electronegative complex plasma. The homotopy analysis method 
(HAM) was developed by Liao \cite{Liao} in 1992 and, applied by Jafari and Seifi \cite{Jafari1} to 
solve the linear and nonlinear fractional diffusion-wave equation. The importance of the HAM lies 
in the fact that it admits ADM and HPM as particular cases.

Time-fractional Navier-Stokes equations have been widely studied.
These equations model the motion of a fluid described by many physical phenomena, for example, the blood flow, ocean current, the flow of liquid in pipes, and airflow around the arms of an aircraft \cite{Prakash}.
El-Shahed and Salem \cite{Shahed} generalized the classical Navier-Stokes equations by replacing
the first time derivative by a Caputo fractional derivative of order $\alpha$, where $0<\alpha\leq{1}$.
The authors obtained the exact solution for three different types of equations using
Hankel transform, Fourier sine transform and Laplace transform. Momani and Odibat \cite{Momani}
solved a time-fractional Navier-Stokes equation applying ADM. Ganji et al. \cite{Ganji}
used analytical technique, HPM, for solving the time-fractional Navier-Stokes equation in polar 
coordinates, and the solutions obtained were expressed in a closed-form. Singh and Kumar \cite{Singh} 
adopted the fractional reduced differential transformation method (FRDTM) to obtain an approximate 
analytical solution of time-fractional order multi-dimensional Navier-Stokes equation. Jaber and
Ahmad \cite{Jaber} used residual power series (RPS) method to find the solution of the nonlinear
time-fractional Navier-Stokes equation in two dimensions. Zhang and Wang \cite{Zhang} proposed 
numerical approximation for a class of Navier-Stokes equations with time fractional derivatives.

Some authors combine two powerful methods to obtain another solution method to solve 
equations and systems time-fractional Navier-Stokes equations. Below we describe some of 
these combinations: Mahmood et al. \cite{Mahmood} used the Laplace Adomian Decomposition 
Method (LADM), a combination of the Laplace transform and ADM; Kumar et al. \cite{Kumar} introduced
homotopy perturbation transform method (HPTM), combined Laplace transform with HPM and,
solved a time-fractional Navier-Stokes equation in a tube. Jena and Chakraverty 
\cite{Jena} applied the homotopy perturbation Elzaki transform method (HPETM) and this 
method consists in the combination of Elzaki transform method and HPM; Prakash et. al
\cite{Prakash} proposed $q$-homotopy analysis transform method ($q$-HATM) to obtain a solution 
of coupled fractional Navier-Stokes equation. This method combine the Laplace transform
and HAM.

The time-fractional model for Navier-Stokes equations then has the form of the operator equation 
\cite{Ganji, Kumar, Momani}
\begin{eqnarray}
\label{Navier-Stokes}
\left\{
\begin{array}{lcl}
\displaystyle {^{\rm{C}}_{a}\mathds{D}_{t}^{\alpha,\psi}}u+(u\cdot\nabla)u=-\frac{1}{\rho}\nabla{p}+
\nu\nabla^2 u, \quad 0<\alpha<1,\\
\nabla u=0,\\
\end{array}\right.
\end{eqnarray}
where ${^{\rm{C}}_{a}\mathds{D}_{t}^{\alpha,\psi}}$ is the $\psi$-Caputo fractional derivative of order $\alpha$,
$t$ is the time, $u$ is the velocity vector, $p$ is the pressure, $\nu$ is the kinematics viscosity and $\rho$ is 
the density. In this work we consider two special cases. First, we consider unsteady, one-dimensional motion of 
a viscous fluid in a tube. The time-fractional Navier-Stokes equations in cylindrical coordinates
that governs the flow field in the tube \cite{Kumar, Momani} are given by
\begin{eqnarray}
{^{\rm{C}}_{a}\mathds{D}_{t}^{\alpha,\psi}}u=P+\nu\left(\frac{\partial^2 u}{\partial r^2}+
\frac{1}{r}\frac{\partial u}{\partial r}\right), \quad 0<\alpha<1, \label{Navier-Stokes-1}
\end{eqnarray}
subject to the initial condition
\begin{eqnarray}
u(r,0)=f(r),
\end{eqnarray}
where $u=u(r,t)$, $P=-\frac{1}{\rho}\frac{\partial p}{\partial z}$ and $f(r)$ is a function depending
only on $r$. 

In the second case, we discuss a nonlinear system of time-fractional Navier-Stokes equations for 
an incompressible fluid flow \cite{Birajdar, Prakash, Singh} described by
\begin{eqnarray}
\label{Navier-Stokes-2}
\left\{
\begin{array}{lcl}
\displaystyle {^{\rm{C}}_{a}\mathds{D}_{t}^{\alpha,\psi}}u+u\frac{\partial u}{\partial x}+
v\frac{\partial u}{\partial y}={\rho_0}\left(\frac{\partial^2 u}{\partial x^2}+
\frac{\partial^2 u}{\partial y^2}\right)-\frac{1}{\rho}\frac{\partial p}{\partial x},\\
\displaystyle {^{\rm{C}}_{a}\mathds{D}_{t}^{\alpha,\psi}}u+u\frac{\partial v}{\partial x}+v\frac{\partial v}{\partial y}={\rho_0}\left(\frac{\partial^2 v}{\partial x^2}+\frac{\partial^2 v}{\partial y^2}\right)-\frac{1}{\rho}\frac{\partial p}{\partial y},
\end{array}\right.
\quad 0<\alpha<1,
\end{eqnarray}
subject to the initial conditions
\begin{eqnarray}
u(x,y,a)=f(x,y) \quad\quad \mbox{and} \quad\quad v(x,y,a)=h(x,y),
\end{eqnarray}
where $u=u(x,y,t)$, $v=v(x,y,t)$, $\rho, t, p$ denote constant density, time and pressure, respectively. 
$x, y$ are the spatial components, $\eta$ is the dynamic viscosity, $\rho_0=\eta/\rho$ is the kinematic 
viscosity of the flow, $f(x,y)$ and $h(x,y)$  are two functions depending only on $x$ and $y$.

There are many definitions for the fractional derivative \cite{Kilbas, Oliveira, Sales, Vanterler}. In this
work, we consider the $\psi$-Caputo fractional derivative \cite{Almeida} to discuss the time-fractional
Navier-Stokes equations by means of HAM. This fractional derivative admits as particular cases the classical
Caputo fractional derivative and Caputo-Hadamard fractional derivative and the fact that the derivative 
of a constant is identically zero.

This work has been organized as follows: In Section \ref{sec:2} we present notations and properties
associated with fractional calculus that will use in the remainder of the text. In Section \ref{sec:3}, 
the HAM has been described. In Section \ref{sec:4}, the HAM has been used to solve time-fractional
Navier-Stokes equations. Concluding remarks close the paper.
\section{Fractional calculus}
\label{sec:2}

In this section we present the definitions and some properties of the fractional
integrals and fractional derivatives of a function $f$ with respect to another function $\psi$.
Some of these definitions and properties can be found in \cite{Almeida} and 
\cite{Kilbas}.
\begin{definition}
Let $\alpha>0$, $I=[a,b]$ be a finite or infinite interval, $f$ an integrable function defined on $I$ and 
$\psi\in{C}^1(I)$ an increasing function such that $\psi'(x)\neq{0}$, for all $x\in I$. The left fractional
integral of $f$ with respect to another function $\psi$ of order $\alpha$ is defined as
\textnormal{\cite{ Almeida, Kilbas}}
\begin{eqnarray}
{_{a}{\mathds{I}}^{\alpha,\psi}_t}f(x,t)]=\frac{1}{\Gamma(\alpha)}\int_{a}^{t}\psi'(\tau)(\psi(t)-\psi(\tau))^{\alpha-1}f(x,\tau)\textnormal{d}\tau. \label{integral}
\end{eqnarray}
For $\alpha=0$, we have
$${_{a}{\mathds{I}}^{0,\psi}_t}[f(x,t)]=f(x,t).$$
\end{definition}
\begin{definition}
Let $\alpha>0$, $n\in\mathbb{N}$, $I$ is the interval $-\infty\leq{a}<b\leq{\infty}$, $f,\psi\in{C^n}(I)$ two 
functions such that $\psi$ is increasing and $\psi'(x)\neq{0}$, for all $x\in I$. The left $\psi$-Caputo 
fractional derivative of $f$ of order $\alpha$ is given by \textnormal{\cite{Almeida}}
$${^{\rm{C}}_{a}\mathds{D}_{t}^{\alpha,\psi}}[f(x,t)]={_{a}{\mathds{I}}^{n-\alpha,\psi}_t}
\left(\frac{1}{\psi'(t)}\frac{\partial}{\partial t}\right)^n f(x,t),$$
where
$$n=[\alpha]+1 \quad \mbox{for}\quad \alpha\notin\mathbb{N}, \quad\quad 
n=\alpha\quad \mbox{for} \quad \alpha\in\mathbb{N}.$$
To simplify notation, we will use the abbreviated notation
$$f^{[n],\psi}(x,t)=\left(\frac{1}{\psi'(t)}\frac{\partial}{\partial t}\right)^n f(x,t).$$
\end{definition}
\begin{property}\label{prop-1}
Let $f\in C^{n}[a,b]$, $\alpha>0$ and $\delta>0,$ \textnormal{\cite{Almeida}}.
\begin{enumerate}

\item $f(t)=(\psi(t)-\psi(a))^{\delta-1}$, then

$${_{a}{\mathds{I}}^{\alpha,\psi}_t}f(t)=\frac{\Gamma(\delta)}{\Gamma(\alpha+\delta)}(\psi(t)-\psi(a))^{\alpha+\delta-1}.$$


\item $\displaystyle {_{a}{\mathds{I}}^{\alpha,\psi}_t}{^{\rm{C}}_{a}\mathds{D}_{t}^{\alpha,\psi}}[f(x,t)]=f(x,t)-
\sum_{k=0}^{n-1}\frac{f^{[k],\psi}(x,a)}{k!}(\psi(t)-\psi(a))^k,$ where\\ $n-1<\alpha<n$ with
$n\in\mathbb{N}.$
\end{enumerate}

\end{property}
\begin{definition}
Let $\alpha>0$. The one-parameter Mittag-Leffler function has the power series representation 
\textnormal{\cite{Almeida, Mittag}}
\begin{eqnarray}
E_{\alpha}(t)=\sum_{m=1}^{\infty}\frac{t^m}{\Gamma(m\alpha+1)}, \label{ML}
\end{eqnarray}
where $\Gamma(\cdot)$ is gamma function.
\end{definition}
\section{Homotopy analysis method}
\label{sec:3}

The homotopy analysis method based on the concept of homotopy was proposed by Liao \cite{Liao} and 
the basic idea is to obtain, through an initial guess, an exact solution for linear and nonlinear 
differential equations. In this section, we describe the HAM.

We consider the following nonlinear differential equation in a general form
\begin{eqnarray}
\mathcal{N}[u(x,t)]=0, \label{eq-nl}
\end{eqnarray}
where $\mathcal{N}$ is a nonlinear differential operator, $x$ and $t$ are independent variables and $u$ is
an unknown function. We then construct the so-called zero-order deformation equation
\begin{eqnarray}
(1-p)\mathcal{L}[\varphi(x,t;p)-u_0(x,t)]=phH(x,t)\mathcal{N}[\phi(x,t;p)], \label{ordem-zero}
\end{eqnarray}
where $p\in[0,1]$ is an embedding parameter, $\hbar\neq{0}$ is an auxiliary parameter, $H(x,t)$ is an
auxiliary function and $\phi(x,t;p)$ is a function of $x$, $t$ and $p$. Let $u_0(x,t)$ be an 
initial approximation of Eq.(\ref{eq-nl}) and $\mathcal{L}={^{\rm{C}}_{a}\mathds{D}_{t}^{\alpha,\psi}}$ 
denotes an auxiliary linear differential operator with the property
$$\mathcal{L}[\phi(x,t)]=0, \quad\quad {\mbox{for}} \quad\quad \phi(x,t)=0.$$
When $p=0$ and $p=1$, we have
$$\phi(x,t;0)=u_0(x,t), \qquad {\mbox{and}} \qquad \phi(x,t;1)=u(x,t),$$
respectively. As the embedding parameter $p$ increases from $0$ to $1$, the solution $\phi(x,t;p)$ 
depends upon the embedding parameter $p$ and varies from the initial guess $u_0(x,t)$ 
to the solution $u(x,t)$. 

Expanding $\phi(x,t;p)$ in a Taylor's series with respect to $p$, we have
\begin{eqnarray}
\phi(x,t;p)=u_0(x,t)+\sum_{m=1}^{\infty}u_m(x,t)p^{m}, \label{serie}
\end{eqnarray}
where
$$u_m(x,t)=\frac{1}{m!}\frac{\partial^m}{\partial p^m}\phi(x,t;p)\biggl|_{p=0}.$$
Assume that the auxiliary parameter $\hbar$, the auxiliary function $H(x,t)$, the initial 
approximation $u_0(x,t)$, and the auxiliary linear operator $\mathcal{L}={^{\rm{C}}_{a}\mathds{D}_{t}^{\alpha,\psi}}$ 
are so properly chosen that the series, Eq.(\ref{serie}), converges at $p=1$. Then, 
the series Eq.(\ref{serie}), at $p=1$, becomes
$$u(x,t)=\phi(x,t;1)=u_m(x,t)=u_0(x,t)+\sum_{m=1}^{\infty}u_m(x,t).$$
Differentiating Eq.(\ref{ordem-zero}), $m$ times with respect to $p$, then setting $p=0$, 
and dividing it by $m!$, we obtain the $m$th-order deformation equation
\begin{eqnarray}
\mathcal{L}[u_m(x,t)-\mathcal{X}_m u_{m-1}(x,t)]=\hbar H(x,t)R_m(\vec{u}_{m-1},x,t), \label{mth}
\end{eqnarray}
with $\vec{u}_n=\{u_0(x,t),u_1(x,t),\ldots,u_n(x,t)\}$ and 
$$R_m(\vec{u}_{m-1},x,t)=\frac{1}{(m-1)!}\frac{\partial^{m-1}}{\partial p^{m-1}}\mathcal{N}[\phi(x,t;p)]\biggl|_{p=0}$$
where we have introduced the notation
\begin{eqnarray}
\mathcal{X}_{m}=\left\{
\begin{array}{lcl}
0, \quad m\leq{1},\\
1, \quad m>1. \label{x_m}
\end{array}\right.
\end{eqnarray}
Operating the fractional integral operator ${_{a}}{\mathds{I}}_{t}^{\alpha,\psi}$, given by 
Eq.(\ref{integral}), on both sides of Eq.(\ref{mth}), we have
\begin{eqnarray}
u_m(x,t)&=&\mathcal{X}_{m}u_{m-1}(x,t)-\mathcal{X}_{m}\sum_{k=0}^{n-1}\frac{u_{m-1}^{[k],\psi}(x,a)}{k!}(\psi(t)-\psi(a))^k
\nonumber\\
&+&\hbar H(x,t){_{a}}{\mathds{I}}_{t}^{\alpha,\psi}[R_m(\vec{u}_{m-1},x,t)], \quad\quad m\geq{1}. \label{sol-mth}
\end{eqnarray}
Thus, we obtain $u_1(x,t), u_2(x,t),\cdots$ by means of Eq.(\ref{sol-mth}). 
So, $M$th-order approximation of $u(x,t)$ is given by
$$
u(x,t)=\sum_{m=0}^{M}u_m(x,t),
$$
and for $M\rightarrow\infty$, we get an accurate approximation of Eq.(\ref{eq-nl}).
\section{Applications: time-fractional Navier-Stokes equations}
\label{sec:4}

In this section, we apply the HAM to solve time-fractional Navier-Stokes equations 
in cylindrical coordinates and a system of time-fractional Navier-Stokes equations
with Cartesian coordinates.

\begin{application}
Let $u=u(r,t)$ and $0<\alpha<1$. Consider the time-fractional Navier-Stokes equation 
in cylindrical coordinates \textnormal{\cite{Bairwa, Momani, Ragab}} given by 
\textnormal{Eq.(\ref{Navier-Stokes-1})}, this is,
\begin{eqnarray}
{^{\rm{C}}_{a}\mathds{D}_{t}^{\alpha,\psi}}u=P+\nu\left(\frac{\partial^2 u}{\partial r^2}+
\frac{1}{r}\frac{\partial u}{\partial r}\right) \label{NS1}
\end{eqnarray}
subject to the initial condition
\begin{eqnarray}
u(r,a)=1-r^2. \label{CI1}
\end{eqnarray}
\end{application}
In order to solve Eq.(\ref{NS1}) by means of HAM, satisfying the initial condition given by
Eq.(\ref{CI1}), it is convenient to choose the initial guess
\begin{eqnarray}
u_0(r,t)=1-r^2 \label{init-approx}
\end{eqnarray}
and the linear differential operator 
\begin{eqnarray*}
\mathcal{L}[\phi(r,t;p)]={^{\rm{C}}_{a}\mathds{D}_{t}^{\alpha,\psi}}[\phi(r,t;p)],
\end{eqnarray*}
satisfying the property
$$\mathcal{L}[c]=0,$$
where $c$ is an arbitrary constant. We define the nonlinear differential operator
\begin{eqnarray}
\mathcal{N}[\phi(r,t;p)]={^{\rm{C}}_{a}\mathds{D}_{t}^{\alpha,\psi}}[\phi(r,t;p)]
-\nu\left(\frac{\partial^2}{\partial r^2}\phi(r,t;p)+\frac{1}{r}\frac{\partial}{\partial r}\phi(r,t;p)\right)
-P(1-\mathcal{X}_m). \label{op-nl}
\end{eqnarray}
Using Eq.(\ref{op-nl}) and the assumption $H(x,t)=1$ we construct the zero-order deformation equation
\begin{eqnarray}
(1-p)\mathcal{L}[\phi(r,t;p)-u_0(r,t)]=p\hbar\mathcal{N}[\phi(r,t;p)].
\end{eqnarray}
Obviously, when $p=0$ and $p=1$, we get
\begin{eqnarray*}
\phi(r,t;0)=u_0(r,t) \quad\quad \mbox{and} \quad\quad \phi(r,t;1)=u(r,t),
\end{eqnarray*}
respectively. So the $m$th-order deformation equation is
\begin{eqnarray}
\mathcal{L}[u_m(r,t)-\mathcal{X}_{m}u_{m-1}(r,t)]=\hbar R_{m}(\vec{u}_{m-1},r,t), \label{ordem-m}
\end{eqnarray}
subject to the initial condition $u_m(r,a)=0$ where $\mathcal{X}_m$ is defined by Eq.(\ref{x_m}) and
$$R_{m}(\vec{u}_{m-1},r,t)={^{\rm{C}}_{a}\mathds{D}_{t}^{\alpha,\psi}}u_{m-1}(r,t)-
\nu\left(\frac{\partial^2}{\partial r^2}u_{m-1}(r,t)+\frac{1}{r}\frac{\partial}{\partial r}u_{m-1}(r,t)\right)
-P(1-\mathcal{X}_m).$$
Now we apply the integral fractional operator ${_{a}\mathds{I}_{t}^{\alpha,\psi}}$ on both sides 
of Eq.(\ref{ordem-m}), we have
\begin{eqnarray*}
&&{_{a}\mathds{I}_{t}^{\alpha,\psi}}{^{\rm{C}}_{a}\mathds{D}_{t}^{\alpha,\psi}}
[u_m(r,t)-\mathcal{X}_{m}u_{m-1}(r,t)]\\
&&=\hbar\,{_{a}\mathds{I}_{t}^{\alpha,\psi}}
\left[{^{\rm{C}}_{a}\mathds{D}_{t}^{\alpha,\psi}}u_{m-1}(r,t)-\nu\left(\frac{\partial^2}{\partial r^2}u_{m-1}(r,t)
+\frac{1}{r}\frac{\partial}{\partial r}u_{m-1}(r,t)\right)-P(1-\mathcal{X}_m)\right],
\end{eqnarray*}
and using \textbf{Property \ref{prop-1}}, we obtain the following solution
\begin{eqnarray*}
u_m(r,t)&-&\sum_{k=0}^{n-1}\frac{u_{m}^{[k],\psi}(r,a)}{k!}
(\psi(t)-\psi(a))^k-\mathcal{X}_{m}u_{m-1}(r,t)\\
&+&\mathcal{X}_{m}\sum_{k=0}^{n-1}\frac{u_{m-1}^{[k],\psi}(r,a)}{k!}
(\psi(t)-\psi(a))^k=\hbar\left\{u_{m-1}(r,t)-\sum_{k=0}^{n-1}\frac{u_{m-1}^{[k],\psi}(r,a)}{k!}
(\psi(t)-\psi(a))^k\right.\\
&-&\left. {_{a}\mathds{I}_{t}^{\alpha,\psi}}\left[\nu\left(\frac{\partial^2}{\partial r^2}u_{m-1}(r,t)+
\frac{1}{r}\frac{\partial}{\partial r}u_{m-1}(r,t)\right)+P(1-\mathcal{X}_m)\right]\right\},
\quad\quad m\geq{1}.
\end{eqnarray*}
For $0<\alpha<1$, we have $n=1$ and we can rewrite the last equation as
\begin{eqnarray}
u_m(x,t)&=&(\mathcal{X}_m+\hbar)u_{m-1}(r,t)-(\mathcal{X}_{m}+\hbar)u_{m-1}(r,a) \label{u_m}\\
&-&\hbar\,{_{a}\mathds{I}_{t}^{\alpha,\psi}}
\left[\nu\left(\frac{\partial^2}{\partial r^2}u_{m-1}(r,t)+
\frac{1}{r}\frac{\partial}{\partial r}u_{m-1}(r,t)\right)+P(1-\mathcal{X}_m)\right], \quad\quad m\geq{1}. \nonumber
\end{eqnarray}
From Eq.(\ref{init-approx}) and Eq.(\ref{u_m}), we obtain
\begin{eqnarray*}
u_0(r,t)&=&1-r^2,\\
u_1(r,t)&=&-\hbar\,{_{a}\mathds{I}_{t}^{\alpha,\psi}}\left[\nu\left(\frac{\partial^2}{\partial r^2}u_{0}(r,t)+
\frac{1}{r}\frac{\partial}{\partial r}u_{0}(r,t)\right)+P\right]\\
&=&-\hbar(P-4\nu)\frac{(\psi(t)-\psi(a))^{\alpha}}{\Gamma(\alpha+1)},\\
u_2(r,t)&=&(1+\hbar)u_1(r,t)-\hbar\,{_{a}\mathds{I}_{t}^{\alpha,\psi}}\left[\nu\left(\frac{\partial^2}{\partial r^2}u_{1}(r,t)+
\frac{1}{r}\frac{\partial}{\partial r}u_{1}(r,t)\right)\right]\\
&=&-(1+\hbar)\hbar(P-4\nu)\frac{(\psi(t)-\psi(a))^{\alpha}}{\Gamma(\alpha+1)},\\
u_3(r,t)&=&(1+\hbar)u_2(r,t)-\hbar\,{_{a}\mathds{I}_{t}^{\alpha,\psi}}\left[\nu\left(\frac{\partial^2}{\partial r^2}u_{2}(r,t)+
\frac{1}{r}\frac{\partial}{\partial r}u_{2}(r,t)\right)\right]\\
&=&-(1+\hbar)^2\hbar(P-4\nu)\frac{(\psi(t)-\psi(a))^{\alpha}}{\Gamma(\alpha+1)},\\
&\vdots&
\end{eqnarray*}
An accurate approximation of Eq.(\ref{NS1}) is given by
\begin{eqnarray*}
u(r,t)&=&u_0(r,t)+u_1(r,t)+u_2(r,t)+u_3(r,t)+\cdots\\
&=&1-r^2-\hbar(P-4\nu)\frac{(\psi(t)-\psi(a))^{\alpha}}{\Gamma(\alpha+1)}\sum_{j=0}^{\infty}(1+\hbar)^j.
\end{eqnarray*}
From geometric series, the above series converges for all $\hbar$ in $|1-\hbar|<1$ and 
can rewrite the above equation as
\begin{eqnarray}
u(r,t)=1-r^2+(P-4\nu)\frac{(\psi(t)-\psi(a))^{\alpha}}{\Gamma(\alpha+1)}, \label{sol}
\end{eqnarray}
which is the exact solution. Observe that the series is independent of $\hbar$ .
There are two important special cases of Eq.(\ref{sol}). First, taking $\psi(t)=t$, $a=0$ and
$\nu=1$. In this case, the solution takes the form
\begin{eqnarray}
u(r,t)=1-r^2+\frac{P-4}{\Gamma(\alpha+1)}\,t^{\alpha}. \label{sol-t} 
\end{eqnarray}
The solution given by Eq.(\ref{sol-t}) is the same found by Momani and Odibat using the
ADM \cite{Momani}, by Ragab et al. using HAM \cite{Ragab} 
and by Bairwa and Singh using iterative Laplace transform \cite{Bairwa}.

On the other hand, if $\psi(t)=\ln t$, $a>0$ and $\nu=1$, the solution given by Eq.(\ref{sol}) becomes
\begin{eqnarray}
u(r,t)=1-r^2+\frac{P-4}{\Gamma(\alpha+1)}\left(\ln\frac{t}{a}\right)^{\alpha}. \label{sol-lnt}
\end{eqnarray}
\begin{figure}[H]
\centering
\includegraphics[width=0.65\textwidth]{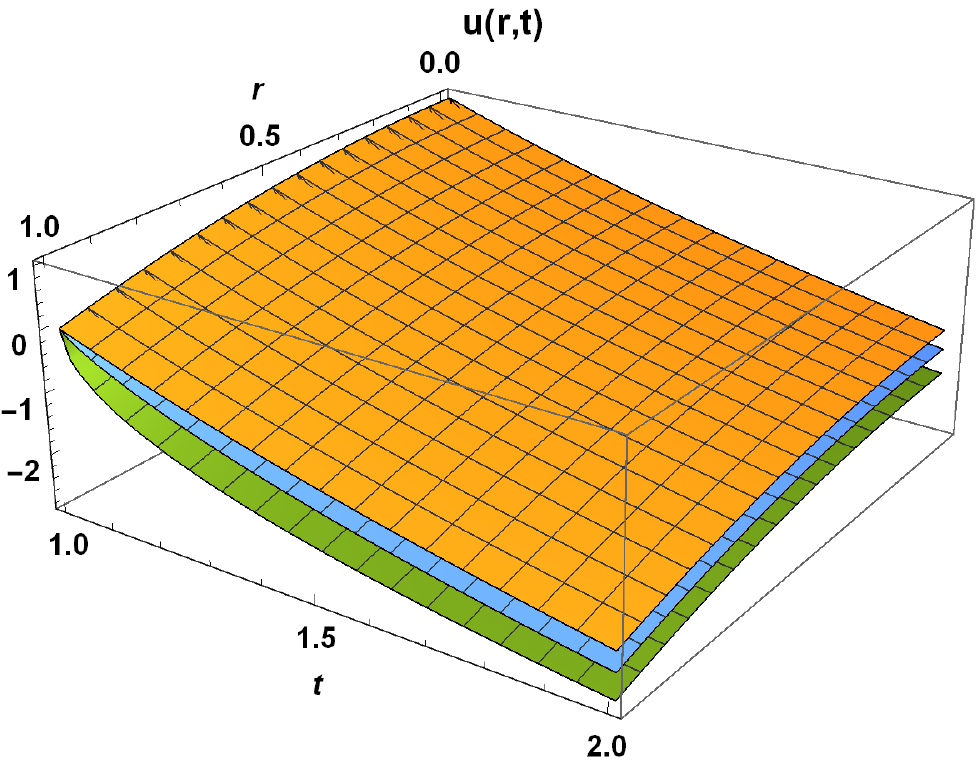}
\caption{Surface plots of the exact solution, Eq.(\ref{sol-lnt}), when $a=P=1$. Orange: $\alpha\rightarrow{1}$;
Blue: $\alpha=0.8$ and; Green: $\alpha=0.5$.}
\label{fig1}
\end{figure}
\begin{figure}[H]
\centering
\includegraphics[width=0.6\textwidth]{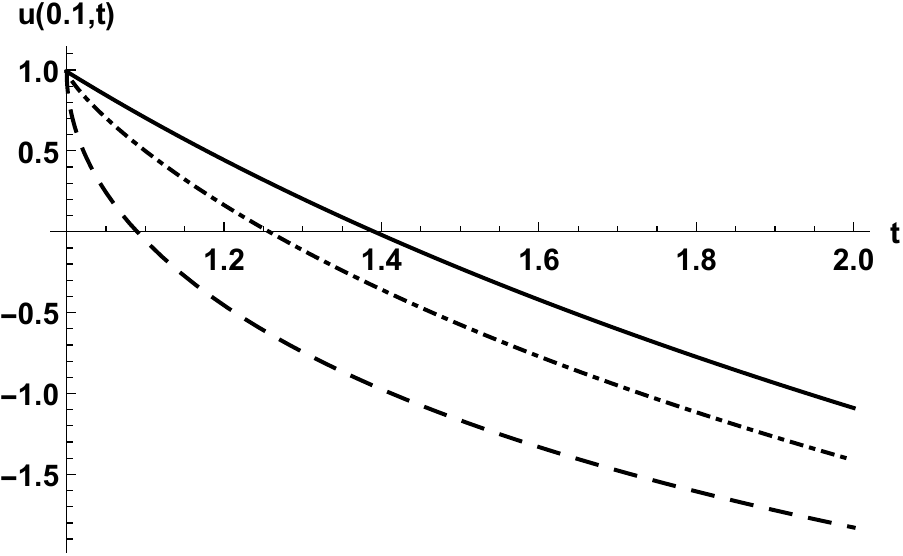}
\caption{Plots of the exact solution, Eq.(\ref{sol-lnt}), when $a=P=1$, $r=0.1$. Solid line: $\alpha\rightarrow{1}$;
Dashdotted: $\alpha=0.8$ and; Dashed: $\alpha=0.5$.}
\label{fig2}
\end{figure}
\begin{application}
Let $u=u(r,t)$ and $0<\alpha<1$. Consider the following time-fractional Navier-Stokes equation 
in cylindrical coordinates \textnormal{\cite{Bairwa, Momani, Ragab}}, given by 
\textnormal{Eq.(\ref{Navier-Stokes})} with $P=0$, this is,
\begin{eqnarray}
{^{\rm{C}}_{a}\mathds{D}_{t}^{\alpha,\psi}}u=\nu\left(\frac{\partial^2 u}{\partial r^2}+
\frac{1}{r}\frac{\partial u}{\partial r}\right) \label{NS2}
\end{eqnarray}
subject to the initial condition
\begin{eqnarray}
u(r,a)=r. \label{CI2}
\end{eqnarray}
\end{application}
In order to solve Eq.(\ref{NS2}) by means of HAM, satisfying the initial condition given by
Eq.(\ref{CI2}), it is convenient to choose the initial guess
\begin{eqnarray}
u_0(r,t)=r. \label{init-approx-2}
\end{eqnarray}
Now we define the nonlinear differential operator by
\begin{eqnarray}
\mathcal{N}[\phi(r,t;p)]={^{\rm{C}}_{a}\mathds{D}_{t}^{\alpha,\psi}}[\phi(r,t;p)]
-\nu\left(\frac{\partial^2}{\partial r^2}\phi(r,t;p)+\frac{1}{r}\frac{\partial}{\partial r}\phi(r,t;p)\right). \label{op-nl-2}
\end{eqnarray}
Using Eq.(\ref{op-nl-2}) and the assumption $H(x,t)=1$ we construct the zero-order deformation equation
\begin{eqnarray}
(1-p)\mathcal{L}[\phi(r,t;p)-u_0(r,t)]=p\hbar\mathcal{N}[\phi(r,t;p)], \label{zero-order-Ap2}
\end{eqnarray}
where the linear differential operator $\mathcal{L}={^{\rm{C}}_{a}\mathds{D}_{t}^{\alpha,\psi}}$
satisfies the property $\mathcal{L}[c]=0$. Obviously, when $p=0$ and $p=1$, the Eq.(\ref{zero-order-Ap2}),
yields
\begin{eqnarray*}
\phi(r,t;0)=u_0(r,t) \quad\quad \mbox{and} \quad\quad \phi(r,t;1)=u(r,t),
\end{eqnarray*}
respectively. So the $m$th-order deformation equation is
\begin{eqnarray}
\mathcal{L}[u_m(r,t)-\mathcal{X}_{m}u_{m-1}(r,t)]=\hbar R_{m}(\vec{u}_{m-1},r,t), \label{ordem-m2}
\end{eqnarray}
subject to the initial condition $u_m(r,a)=0$ where $\mathcal{X}_m$ is defined by Eq.(\ref{x_m}), where
$$R_{m}(\vec{u}_{m-1},r,t)={^{\rm{C}}_{a}\mathds{D}_{t}^{\alpha,\psi}}u_{m-1}-
\nu\left(\frac{\partial^2 u_{m-1}}{\partial r^2}+\frac{1}{r}\frac{\partial u_{m-1}}{\partial r}\right).$$
For $0<\alpha<1$, we have $n=1$ and rearrange some of the terms, thus obtaining
\begin{eqnarray}
u_m(x,t)&=&(\mathcal{X}_m+\hbar)u_{m-1}(r,t)-(\mathcal{X}_{m}+\hbar)u_{m-1}(r,a) \label{u_m2}\\
&-&\hbar\nu\,{_{a}\mathds{I}_{t}^{\alpha,\psi}}
\left[\frac{\partial^2 u_{m-1}}{\partial r^2}+
\frac{1}{r}\frac{\partial u_{m-1}}{\partial r}\right], \quad\quad m\geq{1}. \nonumber
\end{eqnarray}
From Eq.(\ref{init-approx-2}) and Eq.(\ref{u_m2}), we obtain
\begin{eqnarray*}
u_0(r,t)&=&r,\\
u_1(r,t)&=&-\hbar\nu\,{_{a}\mathds{I}_{t}^{\alpha,\psi}}\left[\frac{\partial^2}{\partial r^2}u_{0}(r,t)+
\frac{1}{r}\frac{\partial}{\partial r}u_{0}(r,t)\right]\\
&=&-\frac{\hbar\nu}{r}\frac{(\psi(t)-\psi(a))^{\alpha}}{\Gamma(\alpha+1)},\\
u_2(r,t)&=&(1+\hbar)u_1(r,t)-\hbar\nu\,{_{a}\mathds{I}_{t}^{\alpha,\psi}}\left[\frac{\partial^2}{\partial r^2}u_{1}(r,t)+
\frac{1}{r}\frac{\partial}{\partial r}u_{1}(r,t)\right]\\
&=&-(1+\hbar)\frac{\hbar\nu}{r}\frac{(\psi(t)-\psi(a))^{\alpha}}{\Gamma(\alpha+1)}+
\frac{\hbar^2\nu^2}{r^3}\frac{(\psi(t)-\psi(a))^{2\alpha}}{\Gamma(2\alpha+1)},\\
u_3(r,t)&=&(1+\hbar)u_2(r,t)-\hbar\nu\,{_{a}\mathds{I}_{t}^{\alpha,\psi}}\left[\frac{\partial^2}{\partial r^2}u_{2}(r,t)+
\frac{1}{r}\frac{\partial}{\partial r}u_{2}(r,t)\right]\\
&=&-(1+\hbar)^2\frac{\hbar\nu}{r}\frac{(\psi(t)-\psi(a))^{\alpha}}{\Gamma(\alpha+1)}+
2(1+\hbar)\frac{\hbar^2\nu^2}{r^3}\frac{(\psi(t)-\psi(a))^{2\alpha}}{\Gamma(2\alpha+1)}\\
&-&9\frac{\hbar^3\nu^3}{r^5}\frac{(\psi(t)-\psi(a))^{3\alpha}}{\Gamma(3\alpha+1)},\\
u_4(r,t)&=&(1+\hbar)u_3(r,t)-\hbar\nu\,{_{a}\mathds{I}_{t}^{\alpha,\psi}}\left[\frac{\partial^2}{\partial r^2}u_{3}(r,t)+
\frac{1}{r}\frac{\partial}{\partial r}u_{3}(r,t)\right]\\
&=&-(1+\hbar)^3\frac{\hbar\nu}{r}\frac{(\psi(t)-\psi(a))^{\alpha}}{\Gamma(\alpha+1)}+3(1+\hbar)^2\frac{\hbar^2\nu^2}{r^3}
\frac{(\psi(t)-\psi(a))^{2\alpha}}{\Gamma(2\alpha+1)}\\
&-&27(1+\hbar)\frac{\hbar^3\nu^3}{r^5}\frac{(\psi(t)-\psi(a))^{3\alpha}}{\Gamma(3\alpha+1)}+225\frac{\hbar^4\nu^4}{r^7}\frac{(\psi(t)-\psi(a))^{4\alpha}}{\Gamma(4\alpha+1)},\\
&\vdots&
\end{eqnarray*}
An accurate approximation of Eq.(\ref{NS2}) is given by
\begin{eqnarray*}
u(r,t)&=&u_0(r,t)+u_1(r,t)+u_2(r,t)+u_3(r,t)+\cdots\\
&=&r-\frac{\hbar\nu}{r}\frac{(\psi(t)-\psi(a))^{\alpha}}{\Gamma(\alpha+1)}\,
[1+(1+\hbar)+(1+\hbar)^2+(1+\hbar)^3+\cdots]\\
&+&\frac{\hbar^2\nu^2}{r^3}\frac{(\psi(t)-\psi(a))^{2\alpha}}{\Gamma(2\alpha+1)}\,
[1+2(1+\hbar)+3(1+\hbar)^2+4(1+\hbar)^3+\cdots]\\
&-&9\frac{\hbar^3\nu^3}{r^5}\frac{(\psi(t)-\psi(a))^{3\alpha}}{\Gamma(3\alpha+1)}\,
[1+3(1+\hbar)+6(1+\hbar)^2+10(1+\hbar)^3+\cdots]\\
&+&225\frac{\hbar^4\nu^4}{r^7}\frac{(\psi(t)-\psi(a))^{4\alpha}}{\Gamma(4\alpha+1)}\,
[1+4(1+\hbar)+10(1+\hbar)^2+20(1+\hbar)^3+\cdots]+\cdots.\\
\end{eqnarray*}
Note that, the terms in brackets are geometric series and they are convergent for
$|1+\hbar|<1$. Considering $\nu=1$ in the last equation yields
\begin{eqnarray}
u(r,t)&=&r+\frac{1}{r}\frac{(\psi(t)-\psi(a))^{\alpha}}{\Gamma(\alpha+1)}+
\frac{1}{r^3}\frac{(\psi(t)-\psi(a))^{2\alpha}}{\Gamma(2\alpha+1)}
+\frac{9}{r^5}\frac{(\psi(t)-\psi(a))^{3\alpha}}{\Gamma(3\alpha+1)} \nonumber\\
&+&\frac{225}{r^7}\frac{(\psi(t)-\psi(a))^{4\alpha}}{\Gamma(4\alpha+1)}+\cdots \nonumber\\
&=&r+\sum_{k=1}^{\infty}\frac{1^2\times{3^2}\times\cdots\times{(2k-3)^2}}{r^{2k-1}}
\frac{(\psi(t)-\psi(a))^{k\alpha}}{\Gamma(k\alpha+1)} \label{sol-2}
\end{eqnarray}
which is exact solution. The series solution is independent of $\hbar$.
In practice, two special cases of the solution given by Eq.(\ref{sol-2}) are of particular 
importance because they provide the solution of Eq.(\ref{NS2}) subject to the initial 
condition Eq.(\ref{CI2}) considering Caputo and Caputo-Hadamard fractional derivatives.
First, we consider $\psi(t)=t$ and $a=0$ in Eq.(\ref{sol-2}),
\begin{eqnarray}
u(r,t)=r+\sum_{k=1}^{\infty}\frac{1^2\times{3^2}\times\cdots\times{(2k-3)^2}}{r^{2k-1}}
\frac{t^{k\alpha}}{\Gamma(k\alpha+1)}. \label{sol-t-2.1} 
\end{eqnarray}
Eq.(\ref{sol-t-2.1}) is the solution found by Momani and Odibat using the Adomian
decomposition method \cite{Momani}, also is the solution found by Ragab et al. using
HAM \cite{Ragab} and the solution found by Bairwa and Singh using iterative Laplace 
transform method \cite{Bairwa}. In the second case, we consider $\psi(t)=\ln t$ and 
$a>0$, thus the solution Eq.(\ref{sol-2}) becomes
\begin{eqnarray}
u(r,t)=r+\sum_{k=1}^{\infty}\frac{1^2\times{3^2}\times\cdots\times{(2k-3)^2}}{r^{2k-1}}
\frac{1}{\Gamma(k\alpha+1)}\left(\ln\frac{t}{a}\right)^{k\alpha}. \label{sol-t-2.2} 
\end{eqnarray}
\begin{figure}[H]
\centering
\includegraphics[width=0.75\textwidth]{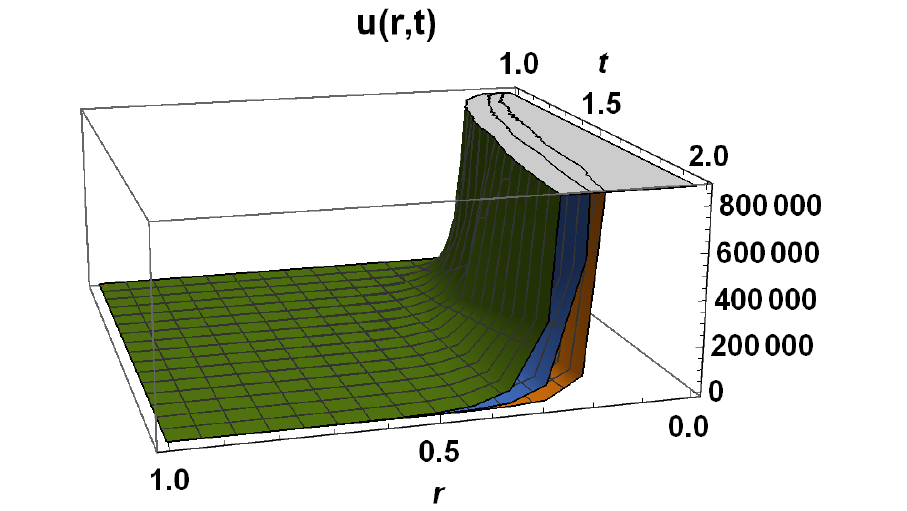}
\caption{Surface plots of approximate solution using 4-terms of Eq.(\ref{sol-t-2.2}) when $a=1$.
Orange: $\alpha\rightarrow{1}$; Blue: $\alpha=0.75$ and; Green: $\alpha=0.5$.}
\label{fig3}
\end{figure}
\begin{figure}[H]
\centering
\includegraphics[width=0.7\textwidth]{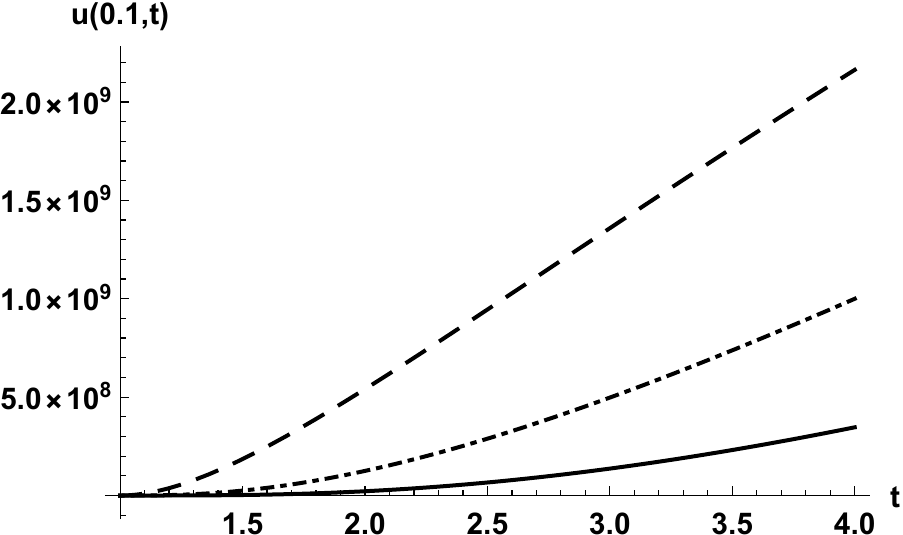}
\caption{Plots of approximate solution using 4-terms of Eq.(\ref{sol-lnt}) when $a=1$ and $r=0.1$.
Solid line: $\alpha\rightarrow{1}$; Dashdotted: $\alpha=0.75$ and; Dashed: $\alpha=0.5$.}
\label{fig4}
\end{figure}
\begin{application}
Let $u=u(x,y,t)$, $v=v(x,y,t)$ and $0<\alpha<1$. Consider a two-dimensional incompressible 
time-fractional Navier-Stokes equations in Cartesian coordinates \textnormal{\cite{Prakash, Singh}}, 
given by \textnormal{Eq.(\ref{Navier-Stokes-2})},this is,
\begin{eqnarray}
\label{NS3}
\left\{
\begin{array}{lcl}
\displaystyle {^{\rm{C}}_{a}\mathds{D}_{t}^{\alpha,\psi}}u+u\frac{\partial u}{\partial x}+v\frac{\partial u}{\partial y}=
{\rho_0}\left(\frac{\partial ^2 u}{\partial x^2}+\frac{\partial^2 u}{\partial y^2}\right)+g,\\
\displaystyle {^{\rm{C}}_{a}\mathds{D}_{t}^{\alpha,\psi}}v+u\frac{\partial v}{\partial x}+v\frac{\partial v}{\partial y}=
{\rho_0}\left(\frac{\partial ^2 v}{\partial x^2}+\frac{\partial^2 v}{\partial y^2}\right)-g, 
\end{array}\right.
\end{eqnarray}
where $g=-\frac{1}{\rho}\frac{\partial p}{\partial x}=\frac{1}{\rho}\frac{\partial p}{\partial y}$ is constant
and subject to the initial conditions
\begin{eqnarray}
u(x,y,a)=-\sin(x+y) \quad\quad \mbox{and} \quad\quad v(x,y,a)=\sin(x+y). \label{CI3}
\end{eqnarray}
\end{application}
In order to solve Eq.(\ref{NS3}) using HAM, satisfying the initial condition given by
Eq.(\ref{CI3}), we choose the following initial guesses
\begin{eqnarray}
u_0(x,y,t)=-\sin(x+y) \quad\quad \mbox{and} \quad\quad v_0(x,y,t)=\sin(x+y) \label{init-approx-3}
\end{eqnarray}
and the linear differential operators
\begin{eqnarray*}
\mathcal{L}_1[\phi_1(x,y,t;p)]&=&{^{\rm{C}}_{a}\mathds{D}_{t}^{\alpha,\psi}}[\phi_i(x,y,t;p)],\\
\mathcal{L}_2[\phi_2(x,y,t;p)]&=&{^{\rm{C}}_{a}\mathds{D}_{t}^{\alpha,\psi}}[\phi_i(x,y,t;p)],
\end{eqnarray*}
satisfying the property $\mathcal{L}[c_i]=0,$ $i=1,2$,
where $c_i$ are arbitrary constants. We define a system of nonlinear operators as
\begin{eqnarray*}
\mathcal{N}_1[\phi_1,\phi_2]&=&{^{\rm{C}}_{a}\mathds{D}_{t}^{\alpha,\psi}}[\phi_1]+{\phi_1}\frac{\partial \phi_1}{\partial x}
+{\phi_2}\frac{\partial \phi_1}{\partial y}-{\rho_0}\left(\frac{\partial^2 \phi_1}{\partial x^2}+
\frac{\partial^2 \phi_1}{\partial y^2}\right)-g,\\
\mathcal{N}_2[\phi_1,\phi_2]&=&{^{\rm{C}}_{a}\mathds{D}_{t}^{\alpha,\psi}}[\phi_2]+{\phi_1}\frac{\partial \phi_2}{\partial x}
+{\phi_2}\frac{\partial \phi_2}{\partial y}-{\rho_0}\left(\frac{\partial^2 \phi_2}{\partial x^2}+
\frac{\partial^2 \phi_2}{\partial y^2}\right)+g.
\end{eqnarray*}
Using Eq.(\ref{op-nl}) and the assumption $H_i(x,t)=1$, $i=1,2$, we construct the zero-order deformation equation
\begin{eqnarray}
(1-p)\mathcal{L}_1[\phi_1(x,y,t;p)-u_0(x,y,t)]=p{\hbar_1}\mathcal{N}_1[\phi_1,\phi_2],\\
(1-p)\mathcal{L}_2[\phi_2(x,y,t;p)-v_0(x,y,t)]=p{\hbar_2}\mathcal{N}_2[\phi_1,\phi_2].
\end{eqnarray}
Obviously, when $p=0$ and $p=1$,
\begin{eqnarray*}
\phi_1(x,y,t;0)=u_0(x,y,t), \quad\quad \phi_1(x,y,t;1)=u(x,y,t),\\
\phi_2(x,y,t;0)=v_0(x,y,t), \quad\quad \phi_2(x,y,t;1)=v(x,y,t).
\end{eqnarray*}
So the $m$th-order deformation equations are given by
\begin{eqnarray}
\mathcal{L}_1[u_m(x,y,t)-\mathcal{X}_{m}u_{m-1}(x,y,t)]=\hbar_1 R_{1m}(\vec{u}_{m-1},\vec{v}_{m-1}), \label{ordem-m3}\\
\mathcal{L}_2[v_m(x,y,t)-\mathcal{X}_{m}v_{m-1}(x,y,t)]=\hbar_2 R_{2m}(\vec{u}_{m-1},\vec{v}_{m-1}),  \label{ordem-m4}
\end{eqnarray}
subject to the initial conditions $u_m(x,y,a)=0$ and $v_m(x,y,a)=0$ where $\mathcal{X}_m$ is defined 
by Eq.(\ref{x_m}) and
\begin{eqnarray*}
&&R_{1m}(\vec{u}_{m-1},\vec{v}_{m-1})={^{\rm{C}}_{a}\mathds{D}_{t}^{\alpha,\psi}}u_{m-1}(x,y,t)+
\sum_{i=0}^{m-1}u_i(x,y,t)\frac{\partial}{\partial x}u_{m-1-i}(x,y,t)\\
&&+\sum_{i=0}^{m-1}v_i(x,y,t)\frac{\partial}{\partial y}u_{m-1-i}(x,y,t)-
\rho_0\left(\frac{\partial^2}{\partial x^2}u_{m-1}(x,y,t)+\frac{\partial}{\partial y^2}u_{m-1}(x,y,t)\right)\\
&&-g(1-\mathcal{X}_m),
\end{eqnarray*}
and
\begin{eqnarray*}
&&R_{2m}(\vec{u}_{m-1},\vec{v}_{m-1})={^{\rm{C}}_{a}\mathds{D}_{t}^{\alpha,\psi}}v_{m-1}(x,y,t)+
\sum_{i=0}^{m-1}u_i(x,y,t)\frac{\partial}{\partial x}v_{m-1-i}(x,y,t)\\
&&+\sum_{i=0}^{m-1}v_i(x,y,t)\frac{\partial}{\partial y}v_{m-1-i}(x,y,t)-
\rho_0\left(\frac{\partial^2}{\partial x^2}v_{m-1}(x,y,t)+\frac{\partial}{\partial y^2}v_{m-1}(x,y,t)\right)\\
&&+g(1-\mathcal{X}_m).
\end{eqnarray*}
Now we apply the integral fractional operator ${_{a}\mathds{I}_{t}^{\alpha,\psi}}$ on both sides 
of Eq.(\ref{ordem-m3}) and Eq.(\ref{ordem-m4}) and using \textbf{Property \ref{prop-1}}, we find
\begin{eqnarray}
&&u_m(x,y,t)=(\mathcal{X}_m+\hbar_1)u_{m-1}(x,y,t)-(\mathcal{X}_{m}+\hbar_1)u_{m-1}(x,y,a) \label{u_m3} \\
&&+\hbar_1\,{_{a}\mathds{I}_{t}^{\alpha,\psi}}
\left[\sum_{i=0}^{m-1}u_i(x,y,t)\frac{\partial}{\partial x}u_{m-1-i}(x,y,t)+
\sum_{i=0}^{m-1}v_i(x,y,t)\frac{\partial}{\partial y}u_{m-1-i}(x,y,t)\right. \nonumber\\
&&-\left.\rho_0\left(\frac{\partial^2}{\partial x^2}u_{m-1}(x,y,t)+\frac{\partial}{\partial y^2}u_{m-1}(x,y,t)\right)
-g(1-\mathcal{X}_m)\right], \quad\quad m\geq{1} \nonumber
\end{eqnarray}
and
\begin{eqnarray}
&&v_m(x,y,t)=(\mathcal{X}_m+\hbar_2)v_{m-1}(x,y,t)-(\mathcal{X}_{m}+\hbar_2)v_{m-1}(x,y,a) \label{u_m4} \\
&&+\hbar_2\,{_{a}\mathds{I}_{t}^{\alpha,\psi}}
\left[\sum_{i=0}^{m-1}u_i(x,y,t)\frac{\partial}{\partial x}v_{m-1-i}(x,y,t)+
\sum_{i=0}^{m-1}v_i(x,y,t)\frac{\partial}{\partial y}v_{m-1-i}(x,y,t)\right. \nonumber\\
&&-\left.\rho_0\left(\frac{\partial^2}{\partial x^2}v_{m-1}(x,y,t)+\frac{\partial}{\partial y^2}v_{m-1}(x,y,t)\right)
+g(1-\mathcal{X}_m)\right], \quad\quad m\geq{1}. \nonumber
\end{eqnarray}
From Eq.(\ref{init-approx-3}), Eq.(\ref{u_m3}) and Eq.(\ref{u_m4}), we obtain
\begin{eqnarray*}
u_0(x,y,t)&=&-\sin(x+y) \quad\quad \mbox{and} \quad\quad v_0(x,y,t)=\sin(x+y),\\
u_1(x,y,t)&=&\hbar_1\,{_{a}\mathds{I}_{t}^{\alpha,\psi}}\left[u_0\frac{\partial u_0}{\partial x}+
v_0\frac{\partial u_0}{\partial y}-{\rho_0}\left(\frac{\partial^2 u_{0}}{\partial x^2}+
\frac{\partial^2 u_{0}}{\partial y^2}\right)-g\right]\\
&=&-2\hbar_1{\rho_0}\sin(x+y)\frac{(\psi(t)-\psi(a))^{\alpha}}{\Gamma(\alpha+1)}-
\hbar_1{g}\frac{(\psi(t)-\psi(a))^{\alpha}}{\Gamma(\alpha+1)},\\
v_1(x,y,t)&=&\hbar_2\,{_{a}\mathds{I}_{t}^{\alpha,\psi}}\left[u_0\frac{\partial v_0}{\partial x}+
v_0\frac{\partial v_0}{\partial y}-{\rho_0}\left(\frac{\partial^2 v_{0}}{\partial x^2}+
\frac{\partial^2 v_{0}}{\partial y^2}\right)+g\right]\\
&=&2\hbar_2{\rho_0}\sin(x+y)\frac{(\psi(t)-\psi(a))^{\alpha}}{\Gamma(\alpha+1)}+
\hbar_2{g}\frac{(\psi(t)-\psi(a))^{\alpha}}{\Gamma(\alpha+1)},\\
u_2(x,y,t)&=&(1+\hbar_1)u_1(x,y,t)+\hbar_1\,{_{a}\mathds{I}_{t}^{\alpha,\psi}}
\left[u_0\frac{\partial u_1}{\partial x}+v_0\frac{\partial u_1}{\partial y}
-{\rho_0}\left(\frac{\partial^2 u_{1}}{\partial x^2}+\frac{\partial^2 u_{1}}{\partial y^2}\right)\right]\\
&=&-2(1+\hbar_1)\hbar_1{\rho_0}\sin(x+y)\frac{(\psi(t)-\psi(a))^{\alpha}}{\Gamma(\alpha+1)}
-(2{\hbar_1\rho_0})^2\sin(x+y)\frac{(\psi(t)-\psi(a))^{2\alpha}}{\Gamma(2\alpha+1)}\\
&-&(1+\hbar_1)\hbar_1\,{g}\frac{(\psi(t)-\psi(a))^{\alpha}}{\Gamma(\alpha+1)},\\
v_2(x,y,t)&=&(1+\hbar_2)v_2(x,y,t)+\hbar_2\,{_{a}\mathds{I}_{t}^{\alpha,\psi}}
\left[u_0\frac{\partial v_1}{\partial x}+v_0\frac{\partial v_1}{\partial y}
-{\rho_0}\left(\frac{\partial^2 v_{1}}{\partial x^2}+\frac{\partial^2 v_{1}}{\partial y^2}\right)\right]\\
&=&-(1+\hbar_2)^2\hbar_2(P-4\nu)\frac{(\psi(t)-\psi(a))^{\alpha}}{\Gamma(\alpha+1)},\\
&=&2(1+\hbar_2)\hbar_2{\rho_0}\sin(x+y)\frac{(\psi(t)-\psi(a))^{\alpha}}{\Gamma(\alpha+1)}
+(2{\hbar_2\rho_0})^2\sin(x+y)\frac{(\psi(t)-\psi(a))^{2\alpha}}{\Gamma(2\alpha+1)}\\
&+&(1+\hbar_2)\hbar_2\,{g}\frac{(\psi(t)-\psi(a))^{\alpha}}{\Gamma(\alpha+1)},\\
u_3(x,y,t)&=&(1+\hbar_1)u_2(x,y,t)+\hbar_1\,{_{a}\mathds{I}_{t}^{\alpha,\psi}}
\left[u_0\frac{\partial u_2}{\partial x}+v_0\frac{\partial u_2}{\partial y}
-{\rho_0}\left(\frac{\partial^2 u_{2}}{\partial x^2}+\frac{\partial^2 u_{2}}{\partial y^2}\right)\right]\\
&=&-2(1+\hbar_1)^2\hbar_1{\rho_0}\sin(x+y)\frac{(\psi(t)-\psi(a))^{\alpha}}{\Gamma(\alpha+1)}
-8(1+\hbar_1)({\hbar_1\rho_0})^2\sin(x+y)\frac{(\psi(t)-\psi(a))^{2\alpha}}{\Gamma(2\alpha+1)}\\
&-&(2{\hbar_1\rho_0})^3\sin(x+y)\frac{(\psi(t)-\psi(a))^{3\alpha}}{\Gamma(3\alpha+1)}-
(1+\hbar_1)^2\hbar_1\,{g}\frac{(\psi(t)-\psi(a))^{\alpha}}{\Gamma(\alpha+1)},\\
v_3(x,y,t)&=&(1+\hbar_2)v_2(x,y,t)+\hbar_2\,{_{a}\mathds{I}_{t}^{\alpha,\psi}}
\left[u_0\frac{\partial v_2}{\partial x}+v_0\frac{\partial v_2}{\partial y}
-{\rho_0}\left(\frac{\partial^2 v_{2}}{\partial x^2}+\frac{\partial^2 v_{2}}{\partial y^2}\right)\right]\\
&=&2(1+\hbar_2)^2\hbar_2{\rho_0}\sin(x+y)\frac{(\psi(t)-\psi(a))^{\alpha}}{\Gamma(\alpha+1)}
+8(1+\hbar_2)({\hbar_2\rho_0})^2\sin(x+y)\frac{(\psi(t)-\psi(a))^{2\alpha}}{\Gamma(2\alpha+1)}\\
&+&(2{\hbar_2\rho_0})^3\sin(x+y)\frac{(\psi(t)-\psi(a))^{3\alpha}}{\Gamma(3\alpha+1)}+
(1+\hbar_2)^2\hbar_2\,{g}\frac{(\psi(t)-\psi(a))^{\alpha}}{\Gamma(\alpha+1)},\\
&\vdots&
\end{eqnarray*}
An accurate approximation of Eq.(\ref{NS1}) is given by
\begin{eqnarray*}
u(x,y,t)&=&u_0(x,y,t)+u_1(x,y,t)+u_2(r,t)+u_3(x,y,t)+\cdots\\
&=&-\sin(x+y)\left\{1+2{\rho_0}\hbar_1\,\frac{(\psi(t)-\psi(a))^{\alpha}}{\Gamma(\alpha+1)}\,
[1+(1+\hbar_1)+(1+\hbar_1)^2+(1+\hbar_1^3)+\cdots]\right.\\
&+&(2{\rho_0}\hbar_1)^2\frac{(\psi(t)-\psi(a))^{2\alpha}}{\Gamma(2\alpha+1)}\,
[1+2(1+\hbar_1)+3(1+\hbar_1)^2+4(1+\hbar_1)^3+\cdots]\\
&+&\left. (2{\rho_0}\hbar_1)^3\frac{(\psi(t)-\psi(a))^{3\alpha}}{\Gamma(3\alpha+1)}\,
[1+3(1+\hbar_1)+6(1+\hbar_1)^2+10(1+\hbar_1)^3+\cdots]+\cdots\right\}\\
&-&\hbar_1 g\,\frac{(\psi(t)-\psi(a))^{\alpha}}{\Gamma(\alpha+1)}[1+(1+\hbar_1)+(1+\hbar_1)^2+\cdots]
\end{eqnarray*}
and
\begin{eqnarray*}
v(x,y,t)&=&v_0(x,y,t)+v_1(x,y,t)+v_2(r,t)+v_3(x,y,t)+\cdots.\\
&=&\sin(x+y)\left\{1+2{\rho_0}\hbar_2\,\frac{(\psi(t)-\psi(a))^{\alpha}}{\Gamma(\alpha+1)}\,
[1+(1+\hbar_2)+(1+\hbar_2)^2+(1+\hbar_2^3)+\cdots]\right.\\
&+&(2{\rho_0}\hbar_2)^2\frac{(\psi(t)-\psi(a))^{2\alpha}}{\Gamma(2\alpha+1)}\,
[1+2(1+\hbar_2)+3(1+\hbar_2)^2+4(1+\hbar_2)^3+\cdots]\\
&+&\left. (2{\rho_0}\hbar_2)^3\frac{(\psi(t)-\psi(a))^{3\alpha}}{\Gamma(3\alpha+1)}\,
[1+3(1+\hbar_2)+6(1+\hbar_2)^2+10(1+\hbar_2)^3+\cdots]+\cdots\right\}\\
&+&\hbar_2 g\,\frac{(\psi(t)-\psi(a))^{\alpha}}{\Gamma(\alpha+1)}[1+(1+\hbar_2)+(1+\hbar_2)^2+\cdots].
\end{eqnarray*}
The terms in brackets are geometric series and they are convergent for $|1+\hbar_i|<1$ with $i=1,2$. 
Thus, we can write
\begin{eqnarray*}
u(x,y,t)&=&-\sin(x+y)\left\{1-\frac{2{\rho_0}(\psi(t)-\psi(a))^{\alpha}}{\Gamma(\alpha+1)}+
\frac{[2{\rho_0}(\psi(t)-\psi(a))^{\alpha}]^{2}}{\Gamma(2\alpha+1)}\right.\\
&-&\left.\frac{[2{\rho_0}(\psi(t)-\psi(a))^{\alpha}]^{3}}{\Gamma(3\alpha+1)}+\cdots\right\}
+{g}\frac{(\psi(t)-\psi(a))^{\alpha}}{\Gamma(\alpha+1)},\\
v(x,y,t)&=&\sin(x+y)\left\{1-\frac{2{\rho_0}(\psi(t)-\psi(a))^{\alpha}}{\Gamma(\alpha+1)}+
\frac{[2{\rho_0}(\psi(t)-\psi(a))^{\alpha}]^{2}}{\Gamma(2\alpha+1)}\right.\\
&-&\left.\frac{[2{\rho_0}(\psi(t)-\psi(a))^{\alpha}]^{3}}{\Gamma(3\alpha+1)}+\cdots\right\}
-{g}\frac{(\psi(t)-\psi(a))^{\alpha}}{\Gamma(\alpha+1)}.
\end{eqnarray*}
Infinite sums can be written in terms of the Mittag-Leffler function, Eq.(\ref{ML}), this is,
\begin{eqnarray}
u(x,y,t)&=&-\sin(x+y)\mathds{E}_{\alpha}[-2\rho_0(\psi(t)-\psi(a))^{\alpha}]+
{g}\frac{(\psi(t)-\psi(a))^{\alpha}}{\Gamma(\alpha+1)}, \label{sol-3} \\
v(x,y,t)&=&\sin(x+y)\mathds{E}_{\alpha}[-2\rho_0(\psi(t)-\psi(a))^{\alpha}]-
{g}\frac{(\psi(t)-\psi(a))^{\alpha}}{\Gamma(\alpha+1)}. \label{sol-4}
\end{eqnarray}
Again, we consider two special cases: First, when $\psi(t)=t$ and $a=g=0$ in Eq.(\ref{sol-3})
and Eq.(\ref{sol-4}). In this case, we have
\begin{eqnarray*}
u(x,y,t)&=&-\sin(x+y)\mathds{E}_{\alpha}[-2\rho_0{t}^{\alpha}],\\
v(x,y,t)&=&\sin(x+y)\mathds{E}_{\alpha}[-2\rho_0{t}^{\alpha}],
\end{eqnarray*}
where $\mathds{E}_{\alpha}(\cdot)$ is the Mittag-Leffler function given by Eq.(\ref{ML}).
These solutions are in agreement with the solutions found by Singh and Kumar using FRDTM \cite{Singh}
and also are in agreement with the solutions found by Prakash et al. applying the $q$-HATM \cite{Prakash}.
In the second case, we consider $\psi(t)=\ln t$, $a=1$ and $g=0$ in Eq.(\ref{sol-3})
and Eq.(\ref{sol-4}),
\begin{eqnarray}
u(x,y,t)&=&-\sin(x+y)\mathds{E}_{\alpha}[-2\rho_0({\ln t})^{\alpha}],\label{sol-ln t-3}\\
v(x,y,t)&=&\sin(x+y)\mathds{E}_{\alpha}[-2\rho_0({\ln t})^{\alpha}]. \label{sol-ln t-4}
\end{eqnarray}
\begin{figure}[H]
\centering
\includegraphics[width=0.49\textwidth]{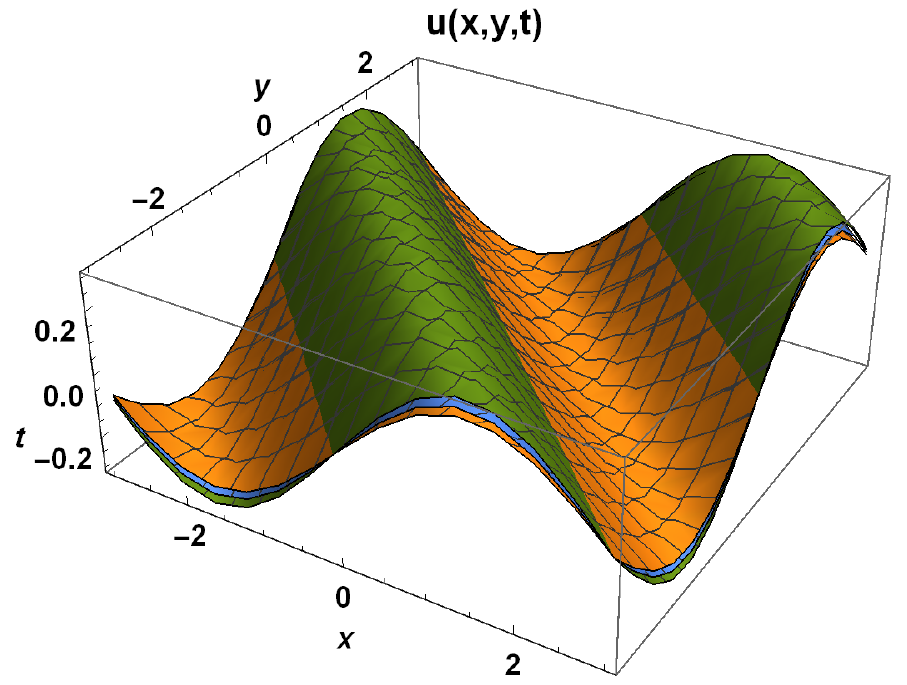}
\includegraphics[width=0.49\textwidth]{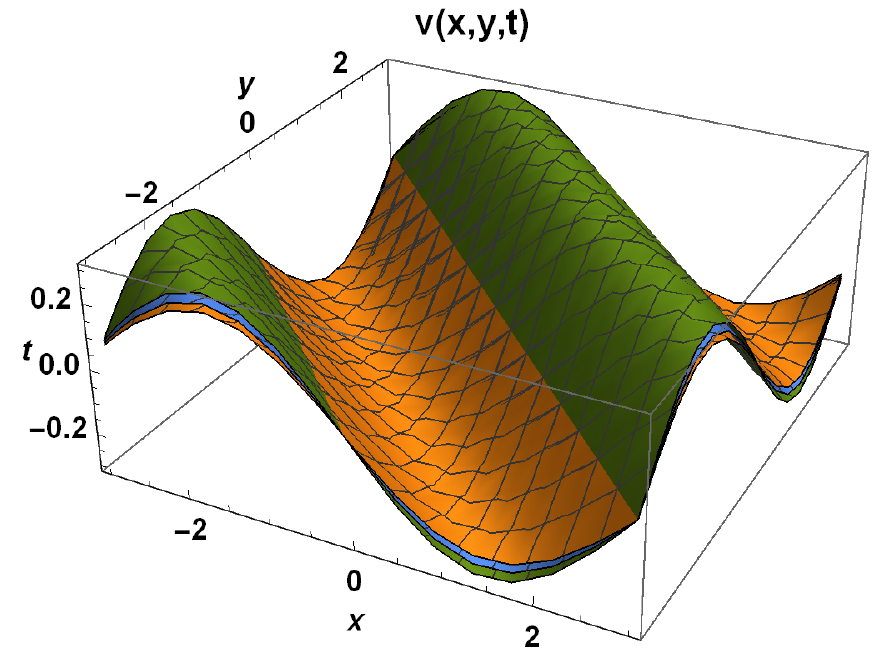}
\caption{Surface plots of exact solutions, Eq.(\ref{sol-ln t-3}) and Eq.(\ref{sol-ln t-4}), for 
$a=1$, $\rho_0=1$ and $t=2$. Orange: $\alpha\rightarrow{1}$; Blue: $\alpha=0.7$ and; Green: $\alpha=0.4$.}
\label{fig5}
\end{figure}
\begin{figure}[H]
\centering
\includegraphics[width=0.49\textwidth]{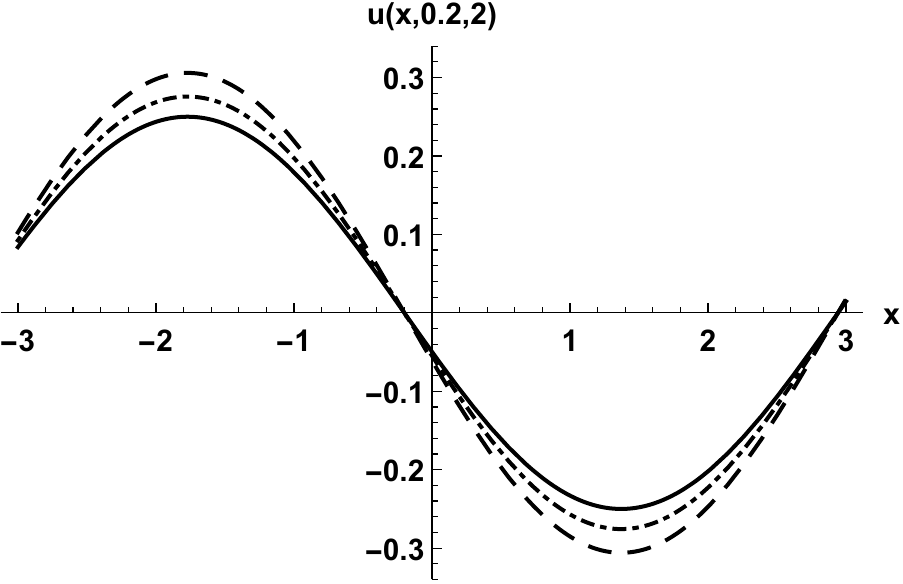}
\includegraphics[width=0.49\textwidth]{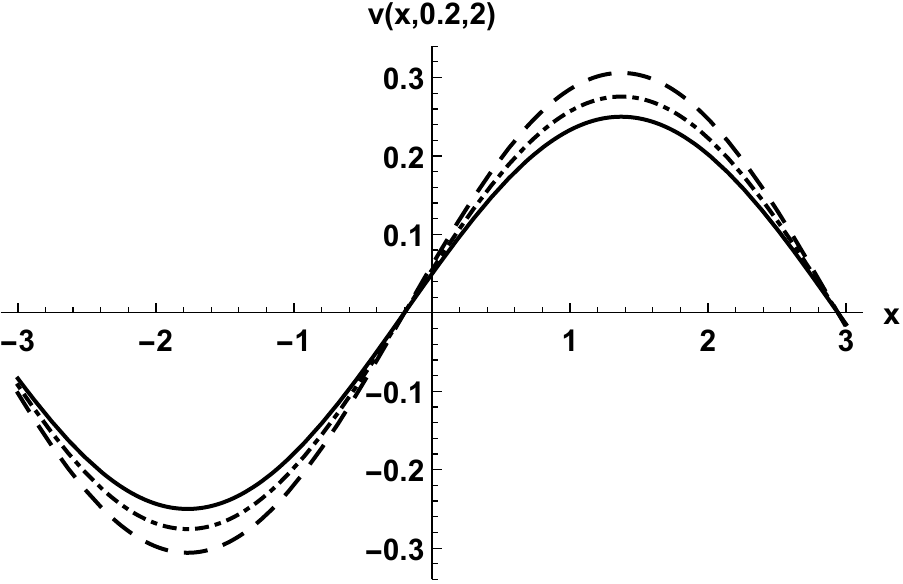}
\caption{Plots exact solution using Eq.(\ref{sol-ln t-3}) and Eq.(\ref{sol-ln t-4}) 
with $a=1$, $\rho_0=1$, $y=0.2$ and $t=2$. Solid line: $\alpha\rightarrow{1}$; Dashdotted: $\alpha=0.7$ 
and; Dashed: $\alpha=0.4$.}
\label{fig6}
\end{figure}
\section{Conclusions}
\label{sec:5}

In this work, HAM was used to solve the unsteady, one-dimensional motion of a viscous 
fluid in a tube which is governed by time-fractional Navier-Stokes equations in cylindrical 
coordinates and, also we solved a nonlinear system of time-fractional Navier-Stokes equations 
for incompressible fluid flow, with Cartesian coordinates. To solve these problems we chose 
the $\psi$-Caputo fractional derivative on time and this operator admits as particular cases 
the Caputo and Caputo-Hadamard fractional derivatives. We graphically represent the solutions 
to these problems when the Caputo-Hadamard fractional derivative is considered. Mathematica 
has been used to draw graphs.



\begin{thebibliography}{}


\bibitem{Adomian}
G. Adomian, Solving Frontier Problems of Physics: The Decomposition Method, Kluwer Acad. Publ., 
Boston, (1994).

\bibitem{Almeida}
R. Almeida, A Caputo fractional derivative of a function with respect to another function, 
Commun. Nonlinear Sci. Numer. Simulat., \textbf{44}:460--481, 2017.

\bibitem{Bairwa}
R. K. Bairwa and J. Singh, Analytical approach to fractional Navier-Stokes equations by iterative Laplace 
transform method, {\textit In}: International workshop of Mathematical Modelling, Applied Analysis and Computation,
J. Singh and D. Kumar, H. Dutta, D. Baleanu and S. Purohit (editors), Springer, \textbf{272}:179--188, (2018).

\bibitem{Birajdar}
G. A. Birajdar, Numerical solution of time fractional Navier-Stokes equation by discrete Adomian decomposition method,
Nonlinear Eng., \textbf{3}:21--26, 2014.

\bibitem{Shahed}
M. El--Shahed and A. Salem, On the generalized Navier-Stokes equations, Appl. Math. Comput.,
\textbf{156}:287--293, 2004.

\bibitem{Ganji}
Z. Z. Ganji, D. D. Ganji, A. D. Ganji and M. Rostamian, Analytical solution of time-fractional Navier-Stokes 
equation in polar coordinate by homotopy perturbation method, Numer. Methods Part. Differ. Equ.,
\textbf{26}:117--124, 2010.

\bibitem{He}
J. H. He, Homotopy perturbation technique, Comput. Math. Appl. Mech. Eng.,
\textbf{178}:257--262, 1999.

\bibitem{Jaber}
K. K. Jaber and R. S. Ahmad, Analytical solution of the time-fractional Navier-Stokes equation,
Ain Shams Eng. J., \textbf{9}:1917--1927, 2018.

\bibitem{Jafari}
H. Jafari and V. Daftardar--Gejji, Solving linear and nonlinear fractional diffusion and wave equations by Adomian decomposition, Appl. Math. Comput., \textbf{180}:488--497", 2006.

\bibitem{Jafari1}
H. Jafari and S. Seifi, Homotopy analysis method for solving linear and nonlinear fractional diffusion-wave equation,
Commun. Nonlinear Sci. Numer. Simulat., \textbf{14}:2006--2012, 2009.

\bibitem{Jena}
R. M. Jena and S. Chakraverty, Solving time-fractional Navier-Stokes equations using homotopy perturbation 
Elzaki transform, SN Appl. Sci., \textbf{1}, 13 pages, 2019.

\bibitem{Kashkari}
B. S. Kashkari, S. A. El--Tantawy, A. H. Salas and L. S. El-Sherif, Homotopy perturbation method for studying dissipative nonplanar solitons in an electronegative complex plasma, Chaos, Solitons \& Fractals, \textbf{130}, 10 pages, 2020.

\bibitem{Kilbas}
A. A. Kilbas, H. M. Srivastava and J. J. Trujillo, Theory and Applications of the Fractional Differential Equations, 
\textbf{204}, Elsevier, Amsterdam, 2006.


\bibitem{Kumar}
D. Kumar, J. Singh and S. Kumar, A fractional model of Navier-Stokes equation arising in unsteady flow of a viscous fluid,
J. Assoc. Arab Univ. Basic Appl. Sci., \textbf{17}:14--19, 2015.

\bibitem{Liao},
S. J. Liao, The proposed homotopy analysis technique for the solution of nonlinear prob-
lems, Ph.D. Thesis, Shanghai Jiao Tong University, (1992).

\bibitem{Mahmood}
S. Mahmood, R. Shah, H. Khan and M. Arif, Laplace Adomian decomposition method for multi dimensional time fractional
model of Navier-Stokes equation, Symmetry, \textbf{11}, 15 pages, 2019.

\bibitem{Mittag}
M. G. Mittag-Leffler, Sur la Nouvelle Fonction ${E}_{\alpha}(x)$, C. R. Acad. Sci., \textbf{137}:554–558, 1903.

\bibitem{Momani}
S. Momani and Z. Odibat, Analytical solution of a time-fractional Navier-Stokes equation by Adomian decomposition method,
Appl. Math. Comput., \textbf{177}:488--494, 2006.

\bibitem{Oliveira}
D. S. Oliveira and E. Capelas de Oliveira, Hilfer-Katugampola fractional derivatives,
Com. Appl. Math., \textbf{37}:3672--3690, 2017.

\bibitem{Prakash}
A. Prakash, P. Veeresha, D. G. Prakasha and M. Goyal, A new efficient technique for solving fractional coupled Navier-Stokes equations using $q$--homotopy analysis transform method, Pramana--J. Phys., \textbf{93}, 10 pages, 2019.

\bibitem{Ragab}
A. A. Ragab, K.M. Hemida, M.S. Mohamed and M.A. Abd El Salam, Solution of time-fractional Navier-Stokes equation by using homotopy analysis method, Gen. Math. Notes, \textbf{13}:13--21, 2012.

\bibitem{Sales}
G. Sales Teodoro, J. A. Tenreiro Machado and E. Capelas de Oliveira, A review of definition of fractional 
derivatives and other operators, J. Comput. Phys., \textbf{388}:195--208, 2019.

\bibitem{Singh}
B. K. Singh and P. Kumar, FRDTM for numerical simulation of multi-dimensional, time-fractional model
of Navier-Stokes equation, Ain Shams Eng. J., \textbf{9}: 827--834, 2018.

\bibitem{Vanterler}
J. Vanterler da C. Sousa and E. Capelas de Oliveira, On the $\psi$-Hilfer fractional derivative,
Commun. Nonlinear Sci. Numer. Simulat., \textbf{60}:72--91, 2018.

\bibitem{Zhang}
J. Zhang and J. Wang, Numerical analysis for Navier-Stokes equations with time fractional derivatives,
Appl Math Comput., \textbf{336}:481--489, 2018.

\end{thebibliography}
\end{document}